\makeatletter
\newcommand{\dontusepackage}[2][]{%
  \@namedef{ver@#2.sty}{9999/12/31}%
  \@namedef{opt@#2.sty}{#1}}
\makeatother
\dontusepackage{subfigure}

\documentclass[]{article}

\usepackage{lmodern}
\usepackage{amssymb,amsmath}
\usepackage{ifxetex,ifluatex}
\usepackage[usenames,dvipsnames]{color}
\usepackage{fixltx2e} % provides \textsubscript
\ifnum 0\ifxetex 1\fi\ifluatex 1\fi=0 % if pdftex
  \usepackage[T1]{fontenc}
  \usepackage[utf8]{inputenc}
\else % if luatex or xelatex
  \ifxetex
    \usepackage{mathspec}
    \usepackage{xltxtra,xunicode}
  \else
    \usepackage{fontspec}
  \fi
  \defaultfontfeatures{Mapping=tex-text,Scale=MatchLowercase}
  
\fi
\IfFileExists{upquote.sty}{\usepackage{upquote}}{}
\IfFileExists{microtype.sty}{%
\usepackage{microtype}
\UseMicrotypeSet[protrusion]{basicmath} % disable protrusion for tt fonts
}{}
\usepackage[margin=1.0in,bottom=1.5in]{geometry}
\usepackage[]{natbib}
\usepackage{listings}
\usepackage{graphicx}
\makeatletter
\def\maxwidth{\ifdim\Gin@nat@width>\linewidth\linewidth\else\Gin@nat@width\fi}
\def\maxheight{\ifdim\Gin@nat@height>\textheight\textheight\else\Gin@nat@height\fi}
\makeatother
\setkeys{Gin}{width=\maxwidth,height=\maxheight,keepaspectratio}
\usepackage{caption}
\usepackage{float}
\usepackage{subfig}
\ifxetex
  \usepackage[setpagesize=false, % page size defined by xetex
              unicode=false, % unicode breaks when used with xetex
              xetex]{hyperref}
\else
  \usepackage[unicode=true]{hyperref}
\fi
\hypersetup{breaklinks=true,
            bookmarks=true,
            pdfauthor={},
            pdftitle={Wavefield recovery with limited-subspace weighted matrix factorizations},
            colorlinks=true,
            citecolor=black,
            urlcolor=blue,
            linkcolor=black,
            pdfborder={0 0 0}}
\def\vector#1{\mathbf{#1}}

\DeclareMathOperator\minim{\hbox{minimize}}

\usepackage{algorithm2e}
\def\minim{\mathop{\hbox{minimize}}}
\title{Wavefield recovery with limited-subspace weighted matrix factorizations}
\author{Yijun Zhang\textsuperscript{1}, Shashin Sharan\textsuperscript{2}, Oscar
Lopez\textsuperscript{4}, Felix J.
Herrmann\textsuperscript{1,2,3}\\\textsuperscript{1} Department of
Electrical \& Computer Engineering, Georgia Institute of
Technology\\\textsuperscript{2} Department of Earth \& Atmospheric
Sciences, Georgia Institute of Technology\\\textsuperscript{3} School of
Computational Science and Engineering, Georgia Institute of
Technology\\\textsuperscript{4} Optimization and Uncertainty
Quantification, Sandia National Laboratories}
\date{}

\begin{document}
\maketitle
\begin{abstract}
Modern-day seismic imaging and monitoring technology increasingly rely
on dense full-azimuth sampling. Unfortunately, the costs of acquiring
densely sampled data rapidly become prohibitive and we need to look for
ways to sparsely collect data, e.g.~from sparsely distributed ocean
bottom nodes, from which we then derive densely sampled surveys through
the method of wavefield reconstruction. Because of their relatively
cheap and simple calculations, wavefield reconstruction via matrix
factorizations has proven to be a viable and scalable alternative to the
more generally used transform-based methods. While this method is
capable of processing all full azimuth data frequency by frequency
slice, its performance degrades at higher frequencies because
monochromatic data at these frequencies is not as well approximated by
low-rank factorizations. We address this problem by proposing a
recursive recovery technique, which involves weighted matrix
factorizations where recovered wavefields at the lower frequencies serve
as prior information for the recovery of the higher frequencies. To
limit the adverse effects of potential overfitting, we propose a
limited-subspace recursively weighted matrix factorization approach
where the size of the row and column subspaces to construct the weight
matrices is constrained. We apply our method to data collected from the
Gulf of Suez, and our results show that our limited-subspace weighted
recovery method significantly improves the recovery quality.
\end{abstract}

\vspace*{-0.45cm}

\section{Introduction}\label{introduction}

Seismic data acquisition plays a key role in the initial phase of oil \&
gas exploration. It also represents a significant budget item for
monitoring of carbon sequestration. For these reasons, it is a challenge
to come up with new acquisition methodologies that improve acquisition
productivity \citep{mosher2014increasing} without sacrificing data
quality. Randomized acquisition according to the principles of
compressive sensing \citep{herrmann2012fighting} in combination with
large-scale wavefield reconstruction algorithms
\citep{kumar2015efficient} has proven a viable tool to improve the
acquisition productivity both in marine and land seismic settings.

So far, many of the employed approached of wavefield reconstruction are
based transform-domain sparsity, which is deigned to explore local
smoothness typically in small windows in up to five dimensions. While
these approaches have been applied successfully on production data, they
do no exploit redundancies present in the data over long distances.
Recovery techniques based on low-rank matrix factorizations
\citep{kumar2015efficient} do not suffer from this shortcoming because
this method works with monochromatic frequency slices that contain data
from the complete survey instead of working within small windows
limiting the apperture. By organizing the data in the appropriate
domain, e.g.~midpoint-offset domain for seismic lines, monochromatic
frequency slices permit approximations in low-rank form, which can be
used to recover fully sample wavefields from subsampled data.

While low-rank factorizations have been employed successfully for low
and midrange frequencies, their performance deteriorates at high
frequencies because monochromatic frequency slices can no longer be
approximated accurately by low-rank factorizations. In this work, we
overcome this problem by using the fact that factorizations at
neighboring frequencies live in close-by subspaces. As described in
early work by \citet{aravkin2013robust}; \citet{eftekhari2018weighted},
this property can be exploited by introducing matrix weights defined in
terms of factorizations of near-by frequency slices. Recent work by
\citet{zhang2019high} took this initial a step further by proposing a
recursive approach where factorizations of frequency slices at lower
frequencies are used as weight for factorizations at the higher
frequencies starting at the low frequencies and working its way up.

While this approach has had some success (see e.g.
\citet{zhang2019high}), there is challenge related to the fact that high
frequencies require higher rank factorizations and this can lead to
overfitting when using this higher rank throughout. We avoid this
overfitting, by adapting the rank of the weighting matrices such that
overfitting is avoided. We do this by actively limiting the row and
column subspaces of the weight matrices. Because we avoid overfitting,
we are able to further improve the wavefield recovery. We also introduce
an alternative formulation where the weight matrices are moved from the
constraint, as in \citet{kumar2015efficient}, to the data misfit
objective, which leads to a significant improvement ($20$ to $25$ times
speedup) computational efficiency.

We organize our paper as follows. First, we review the recursively
weighted wavefield recovery via matrix factorization including the new
formulation where the weight appear in the data misfit term. Next, we
discuss how to limit the subspace of our weighted matrix factorizations.
We conclude by demonstrating our approach on a field data example from
the Gulf of Suez, which shows improved recovery quality compared to
conventional recursively weighted matrix completion.

\vspace*{-0.45cm}

\section{Methodology}\label{methodology}

We start by introducing wavefield reconstruction via weighted matrix
factorization. To improve computational efficiency, we move the weight
matrices to the data misfit term so we no longer have to carry out
numerically expensive weighted projections as in
\citep{aravkin2013robust}. Aside from allowing for a much more
computationally efficient implementation, this alternative formulation
also forms the basis for our limited-subspace approach designed to
prevent overfitting at the low frequencies.

\subsection{Weighted low-rank matrix
factorization}\label{weighted-low-rank-matrix-factorization}

Our proposed extension to wavefield reconstruction via recursively
weighted matrix factorization derives from earlier work by
\citet{kumar2015efficient}, \citet{aravkin2013robust}, and
\citet{zhang2019high}, where we solve
\begin{equation}
\begin{aligned}
& \minim_{\vector{X}_i} \quad \|\vector{Q}\vector{X}_i\vector{W}\|_*\\
& \text{subject to} \quad \|\mathcal{A}(\vector{X}_i) - \vector{b}_i\|_2 \leq \tau
\end{aligned}
\label{eqWlrmf}
\end{equation}
 to within a noise-level dependent data misfit tolerance $\tau$. In this
expression, the matrix $\vector{X}_i$ corresponds to a monochromatic
frequency slice in the midpoint/offset domain (in case of $2$D) at the
$i\mathrm{th}$ frequency ($i \in [1, \cdots , n_f]$ with with $n_f$ the
number of frequencies).

During the wavefield recovery, fully sampled frequency slices are
represented by the complex valued matrix,
$\vector{X} \in \mathbb{C}^{n_f \times n_m \times n_h}$ where $n_m$ is
the number of midpoints and $n_h$ the number of offsets. The symbol
$\mathcal{A}(\cdot)$ stands for the subsampling operator, which collects
monochromatic data at the observed source/receiver combinations into the
vector $\vector{b}_i$. Given these observations, we solve for the fully
sampled $\mathbf{X}_i$ for each frequency by minimizing
equation~\ref{eqWlrmf} with weight matrices $\vector{Q}$ and
$\vector{W}$ given by
\begin{equation}
\vector{Q} = {w}_{1}\vector{U}\vector{U}^H + \vector{U}^\perp \vector{U}^{{\perp}{H}}
\label{eqprojl}
\end{equation}
 and
\begin{equation}
\vector{W} = {w}_{2}\vector{V}\vector{V}^H + \vector{V}^\perp \vector{V}^{{\perp}{H}}.
\label{eqprojr}
\end{equation}
 In these expressions for the weight matrices, the
$\vector{U} \in \mathbb{C}^{n_m \times r}$ and
$\vector{V}\in \mathbb{C}^{n_h \times r}$ are the column and row
subspaces that derive from the low-rank factorization of the nearby
frequency slice. $\vector{U}$ and $\vector{V}$ have orthonormal columns
that span top column and row subspaces of nearby frequency slice.
Because these weight matrices include information on the subspaces of
the current factorization, they serves as prior information aiding the
wavefield recovery via the weighted nuclear norm minimization (denoted
by $\|\mathbf{QXW}\|_\ast=\sum_{j=1}^{r}\sigma_j$ with $\sigma_j$ the
$j^{\text{th}}$ singular value). Depending on whether we have confidence
in the fact that the neighboring frequency slice has an overlapping
subspace, we chose the weights $w_1$ and $w_2$ close to $0$ if we have
confidence and close to $1$ if we do not.

While the above weighted formulation has resulted in major improvements
in the recovery when reliable information on a neighboring frequency
slice is available \citep[\citet{aravkin2013robust}, and
\citet{zhang2019high}]{kumar2015efficient}, the minimization in
equation~\ref{eqWlrmf} is complicated by the presence of the weighting
matrices in the nuclear norm objective. As a result, the minimization
becomes computationally expensive. To avoid this complication, we
replace the optimization variable by
$\vector{\bar{X}}_{i} = \vector{Q}\vector{X}_i\vector{W}$, and rewrite
equation~\ref{eqWlrmf} as
\begin{equation}
\begin{aligned}
& \minim_{\vector{\bar{X}}_i} \quad \|\vector{\bar{X}}_i\|_* \\
& \text{subject to} \quad \|\mathcal{A}({\vector{Q}^{-1}}\vector{\bar{X}}_i{\vector{W}^{-1}}) - \vector{b}_i\|_2 \leq \tau
\end{aligned}
\label{eqWlrmfv2}
\end{equation}
 where the modified weighting matrices
\begin{equation}
{\vector{Q}^{-1}} = \frac{1}{w_1} \vector{U}\vector{U}^H + \vector{U}^\perp \vector{U}^{{\perp}{H}}
\label{eqWrow}
\end{equation}
 and
\begin{equation}
{\vector{W}^{-1}} = \frac{1}{w_2} \vector{V}\vector{V}^H +  \vector{V}^\perp \vector{V}^{{\perp}{H}}
\label{eqWcol}
\end{equation}
 are moved from the objective to the data misfit constraint. To reflect
that we changed the problem, we introduced barred quantities from which
the solution original solution can be readily computed---i.e., we
recover the solution
$\vector{X}_i=\vector{Q}^{-1}\vector{\bar{X}}_i\vector{W}^{-1}$ since
$\vector{\bar{X}}_i=\vector{Q}\vector{X}_i\vector{W}$ solves the above
optimization problem. Compared to equation~\ref{eqWlrmf}, this new
formulation does not require nuclear norm projections onto weighted
matrices while its solution is equivalent to equation~\ref{eqWlrmf}.

Like the original formulation, our new formulation lends also itself to
be cast into a low-rank ($r\ll \max(n_m,n_h)$) factorized form so that
expensive SVDs are avoided in the nuclear norm. After factorization our
wavefield reconstruction involves
\begin{equation}
\begin{aligned}
& \minim_{\vector{\bar L}_i, \vector{\bar R}_i} \quad \frac{1}{2} {\left\| \begin{bmatrix} \vector{\bar L}_i \\ \vector{\bar R}_i \end{bmatrix} \right\|}_F^2 \\
& \text{subject to} \quad \|\mathcal{A}{({\vector{Q}^{-1}} \vector{\bar L}_i \vector{\bar R}_i^{H} {\vector{W}^{-1}})} - \vector{b}_i\|_{2} \leq \epsilon,
\end{aligned}
\label{eqwlrf}
\end{equation}
 where the symbol $^{H}$ denotes the Hermitian transpose and
$\|\cdot\|_F$ is the Frobenius norm (2-norm of the vectorized matrix)
\citep{kumar2015efficient, aravkin2013robust, zhang2019high}. Compared
to the original representation for frequency slices, the above factored
form is compressed since it entails the low-rank pair
$\{\mathbf{\bar L}_i,\,\mathbf{\bar R}_i\}$ ,where
$\vector{\bar{X}}_i=\vector{\bar L}_i \vector{\bar R}_i^{H}$, and does
not rely on storage and manipulation of the original and dense
optimization variable $\mathbf{X}_i$ or $\mathbf{\bar X}_i$. Despite
gains in computation, because of the factored form and redefined data
misfit term, challenges remain with recursive weighted matrix
factorizations \citep{zhang2019high} at the high frequencies and as we
will show these have to do with overfitting.

\subsection{Limited subspace weighted
implementation}\label{limited-subspace-weighted-implementation}

To reduce approximation errors at the high frequencies, we can increase
the rank of the factorization throughout. While increasing the rank
leads to better approximations at the high frequencies adapting this
higher rank at the lower frequencies can lead to overfitting. The
resulting poor reconstructions at the lower frequencies can in turn have
a detrimental effect on the reconstruction at higher frequencies, which
information from the lower frequencies as the recursive algorithm sweeps
from the low to the high frequencies.

By choosing the rank for the limited subspace, we reduce the size of the
subspaces of the weight matrices to prevent overfitting at the lower
frequencies. In equations~\ref{eqprojl}, \ref{eqprojr}, \ref{eqWrow}
and~\ref{eqWcol}, we notice that the size of the weight matrices
$\vector{Q}$ and $\vector{W}$ are independent of rank $r$. Therefore, we
can use a limited subspace to remove the influence of overfitting and
get better results.

By limited subspace, we mean that at a given frequency slice, instead of
using a rank $r$ for row and column subspaces $\vector{U}$ and
$\vector{V}$ respectively, we can use a lower rank $r_s$. In this way,
we can choose higher rank $r$ to reconstruct each frequency but use
lower rank $r_s$ to construct the weight matrices ($\vector{Q}$ and
$\vector{W}$). By choosing smaller rank for the subspaces, we mitigate
the negative influence of overfitting. Therefore, in the
limited-subspace method, we are free to choose smaller values for the
$r_s$ for each frequency slice and higher values for the rank $r$ for
the factorization itself (not for the weights) for each frequency.

\vspace*{-0.45cm}

\section{Numerical Experiments}\label{numerical-experiments}

To demonstrate the advocacy of the proposed method, we use $2$D field
seismic data acquired in the Gulf of Suez with number of sources,
${N}_{s}=355$, and number of receivers, ${N}_{r}=355$. The total number
of time samples in this dataset is ${N}_{t}=1024$ and the sampling
interval is $0.004\, \mathrm{s}$. We use a jittered subsampling
\citep{herrmann2008non} mask to remove $75 \%$ of the sources to obtain
the subsampled data. When data is organized in the midpoint-offset
domain, we know that randomized jittered subsampling method breaks the
inherent low-rank property of seismic data while controlling the largest
gap size of the subsampled data \citep{herrmann2008non}. Controlling
largest gap is important because very large gaps are not suitable for
wavefield reconstruction using sparsity-promotion or low-rank matrix
completion. We use the weighted method as described by
\citet{zhang2019high} to reconstruct frequency slices starting at
$10\, \mathrm{Hz}$ and working our way up to $70\, \mathrm{Hz}$. We use
constant rank across all the frequencies for weight matrices and matrix
factorization. We base these choices for $r_s<r$ on visual inspection of
the recovered frequency slices. To avoid overfitting at lower
frequencies we select rank $r_s$ of the limited subspace constant across
all the frequencies. And to better approximation of higher frequencies
we choose higher rank $r$ across all the frequencies. Combination of
higher rank for matrix factorization and smaller rank for limited
subspace avoid the risk of overfitting and at the same time improves the
data reconstruction quality.

\begin{figure}
\centering
\subfloat[\label{full_data}]{\includegraphics[width=0.500\hsize]{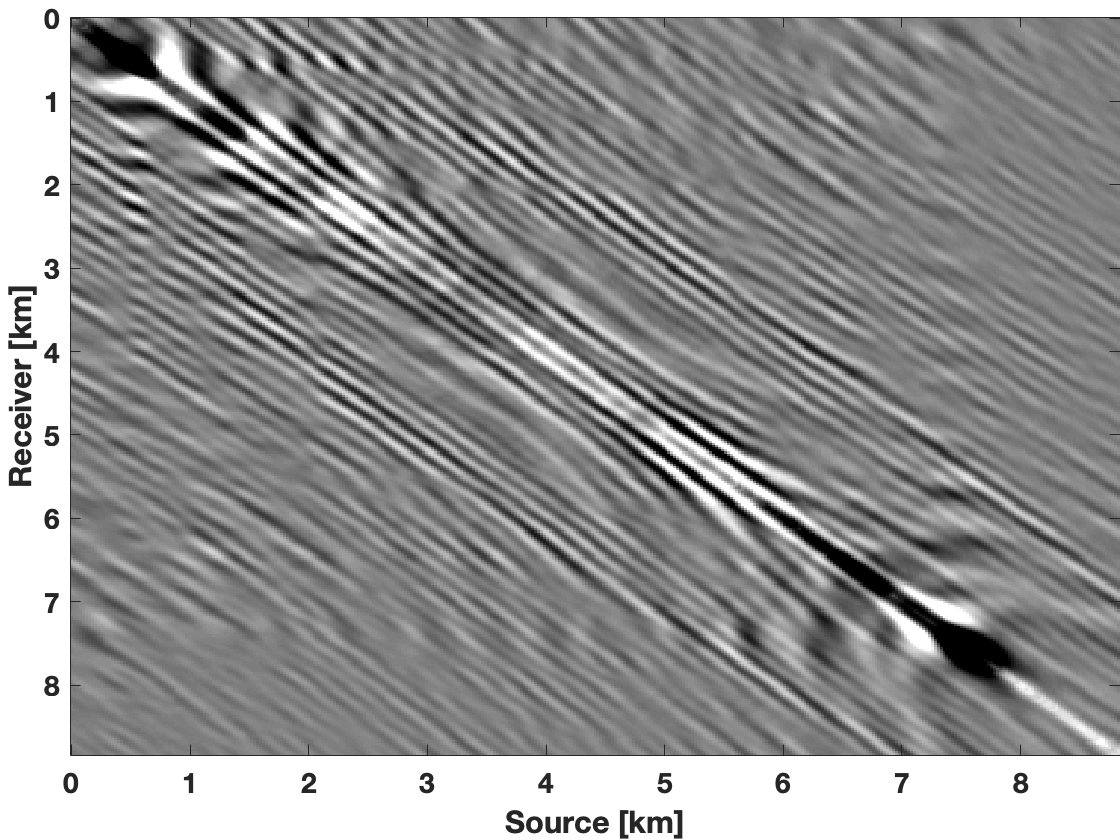}}
\subfloat[\label{jittered_75}]{\includegraphics[width=0.500\hsize]{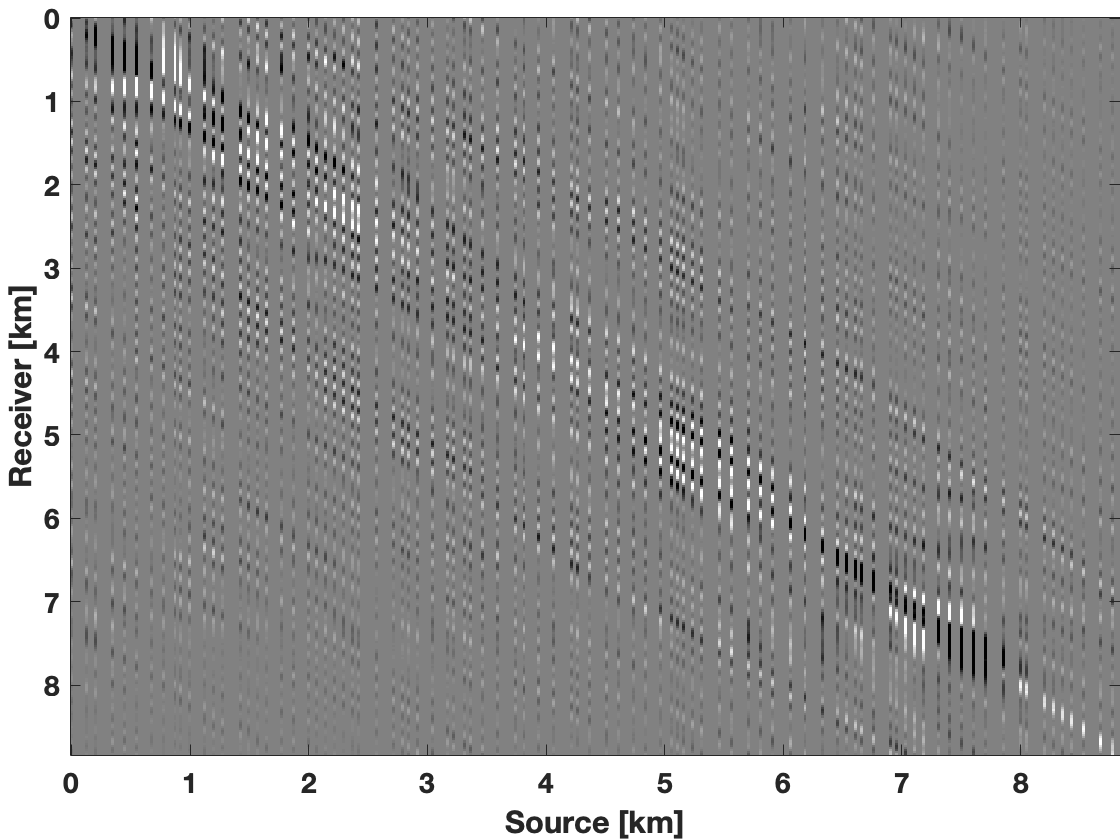}}
\\
\subfloat[\label{full_85}]{\includegraphics[width=0.500\hsize]{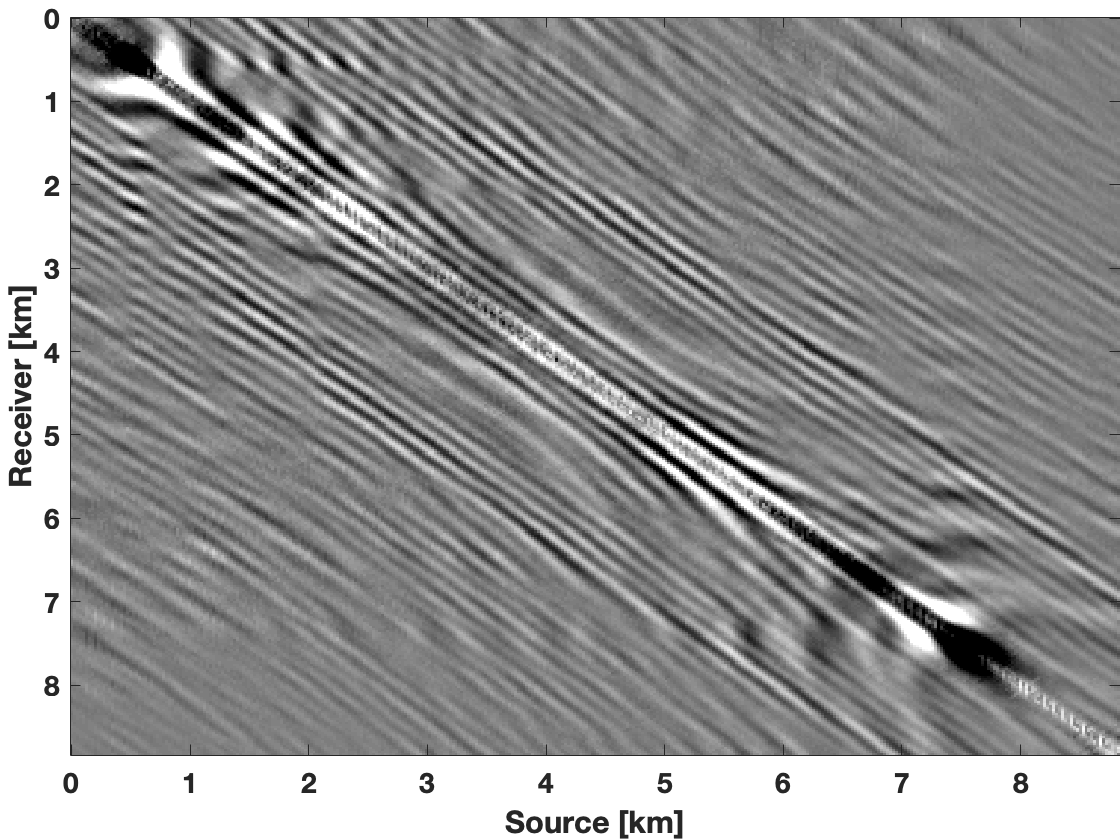}}
\subfloat[\label{full_85_diff}]{\includegraphics[width=0.500\hsize]{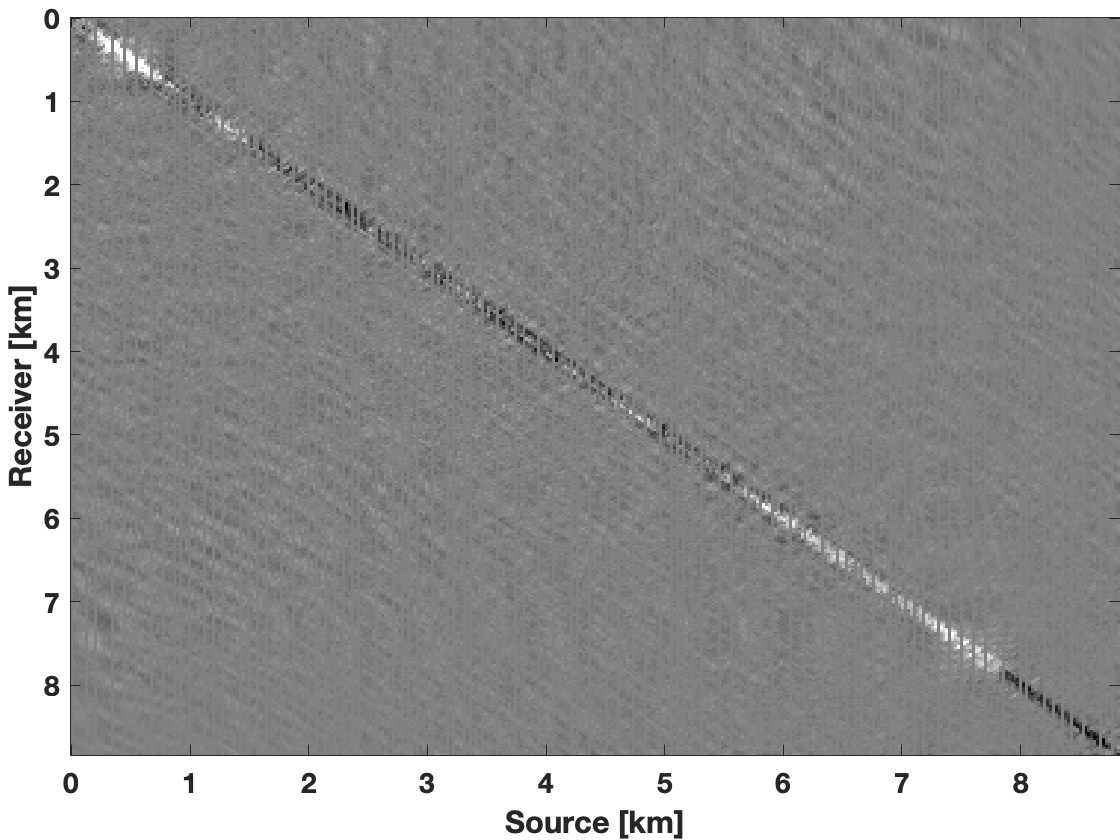}}
\\
\subfloat[\label{full_30}]{\includegraphics[width=0.500\hsize]{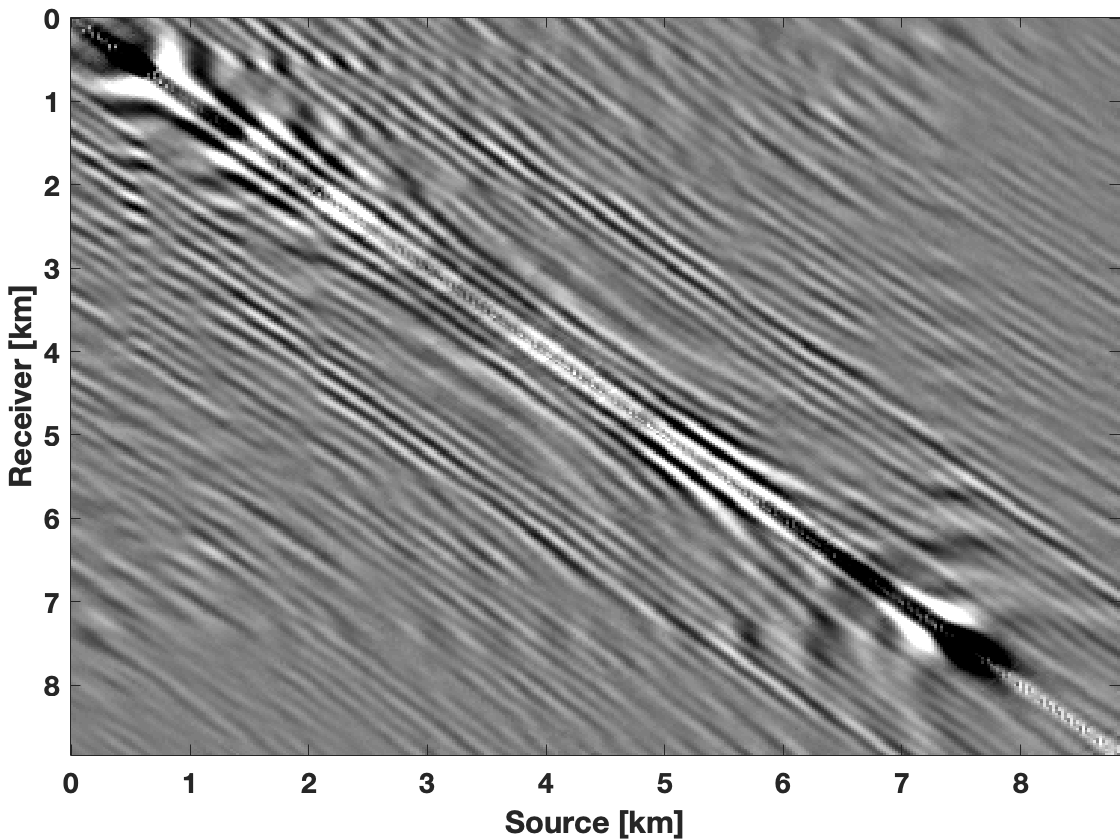}}
\subfloat[\label{full_30_diff}]{\includegraphics[width=0.500\hsize]{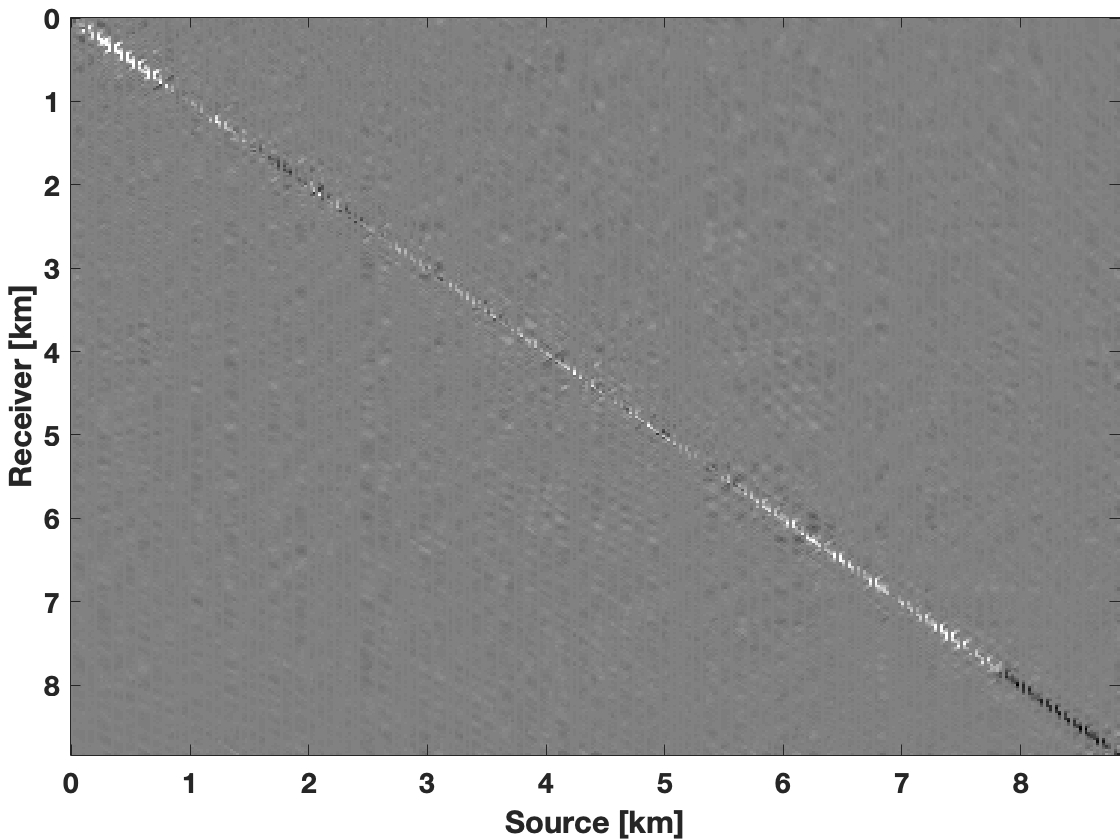}}
\\
\subfloat[\label{sig_SR}]{\includegraphics[width=0.500\hsize]{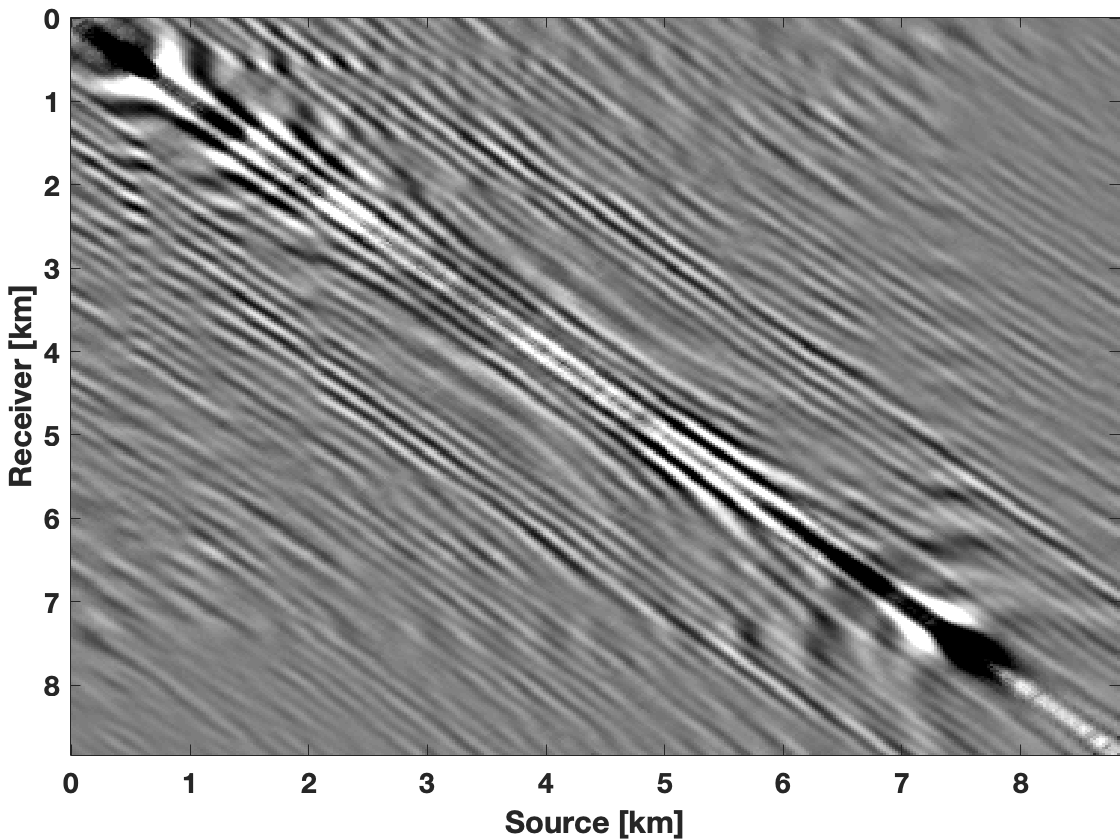}}
\subfloat[\label{sig_diff}]{\includegraphics[width=0.500\hsize]{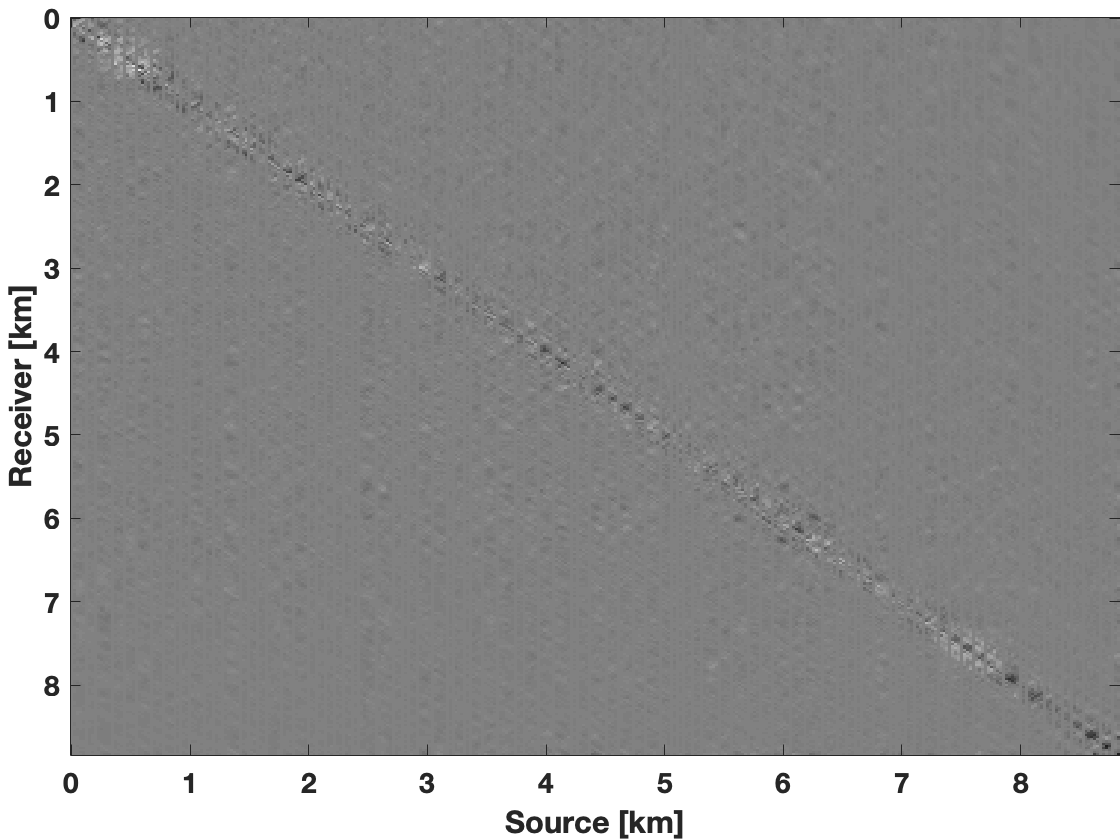}}
\caption{Reconstruction for missing source for a frequency slice at
$22\, \mathrm{Hz}$ shown in the source-receiver domain but reconstructed
in the midpoint-offset domain. \emph{(a)} Ground truth, \emph{(b)}
$75 \%$ subsampled seismic data with jittered subsampling. \emph{(c)}
and \emph{(d)} recovery by weighted matrix factorization
($S/R = 13.09\, \mathrm{dB}$) using conventional recursively weighted
approach with fixed rank $r=85$ and corresponding residual w.r.t. the
ground truth, respectively. \emph{(e)} and \emph{(f)} contain recovery
($S/R = 15.50\, \mathrm{dB}$) for conventional recursively weighted with
a rank $r=25$ and corresponding residual w.r.t. the ground truth
respectively. \emph{(g)} and \emph{(h)} represent recovery
($S/R = 19.52\, \mathrm{dB}$) using limited-subspace weighted method
with limited-subspace rank $r_s = 25$ and corresponding residual w.r.t.
the ground truth respectively.}\label{example_1}
\end{figure}

To demonstrate that the limited-subspace recursively weighted method
gives improved results compared to conventional recursively weighted
method \citep{zhang2019high}, we first show results in the frequency
domain. For each frequency slice, we perform $150$ iterations for both
the methods. For the limited-subspace weighted method, we use rank
$r = 85$ and limited subspace rank of $r_s = 25$. For comparison with
the conventional weighted method, we perform two experiments with a
fixed high rank of $r=85$ and lower rank of $r=25$. We choose lower rank
for conventional weighted method to show that smaller rank itself is not
sufficient for significant improvement in data reconstruction at higher
frequencies. On the other hand we choose higher rank of 85 for
conventional weighted method to show that higher rank is alone not
sufficient to improve the quality of reconstructed data at higher
frequencies because of the overfitting at lower frequencies. We show
reconstruction results for a frequency slice at $22\, \mathrm{Hz}$ in
Figure~\ref{example_1}. Due to overfiting, the conventional method with
rank $r=85$ gives a reconstruction with a smaller S/R of
$13.09\, \mathrm{dB}$ compared to the wavefield reconstruction
(Figures~\ref{full_85} and~\ref{full_85_diff}) obtained with the smaller
rank $r=25$ for which we get S/R of $15.50\, \mathrm{dB}$
(Figures~\ref{full_30} and~\ref{full_30_diff}). We get S/R of
$19.52\, \mathrm{dB}$ for the reconstructed data (Figures~\ref{sig_SR})
using the limited-subspace weighted method. Figure~\ref{sig_diff} shows
the data residual with respect to the ground truth
(Figure~\ref{full_data}). Clearly, our limited-subspace weighted method
outperforms the conventional weighted method in terms of improved
quality of reconstructed data.

\begin{figure}
\centering
\includegraphics[width=1.000\hsize]{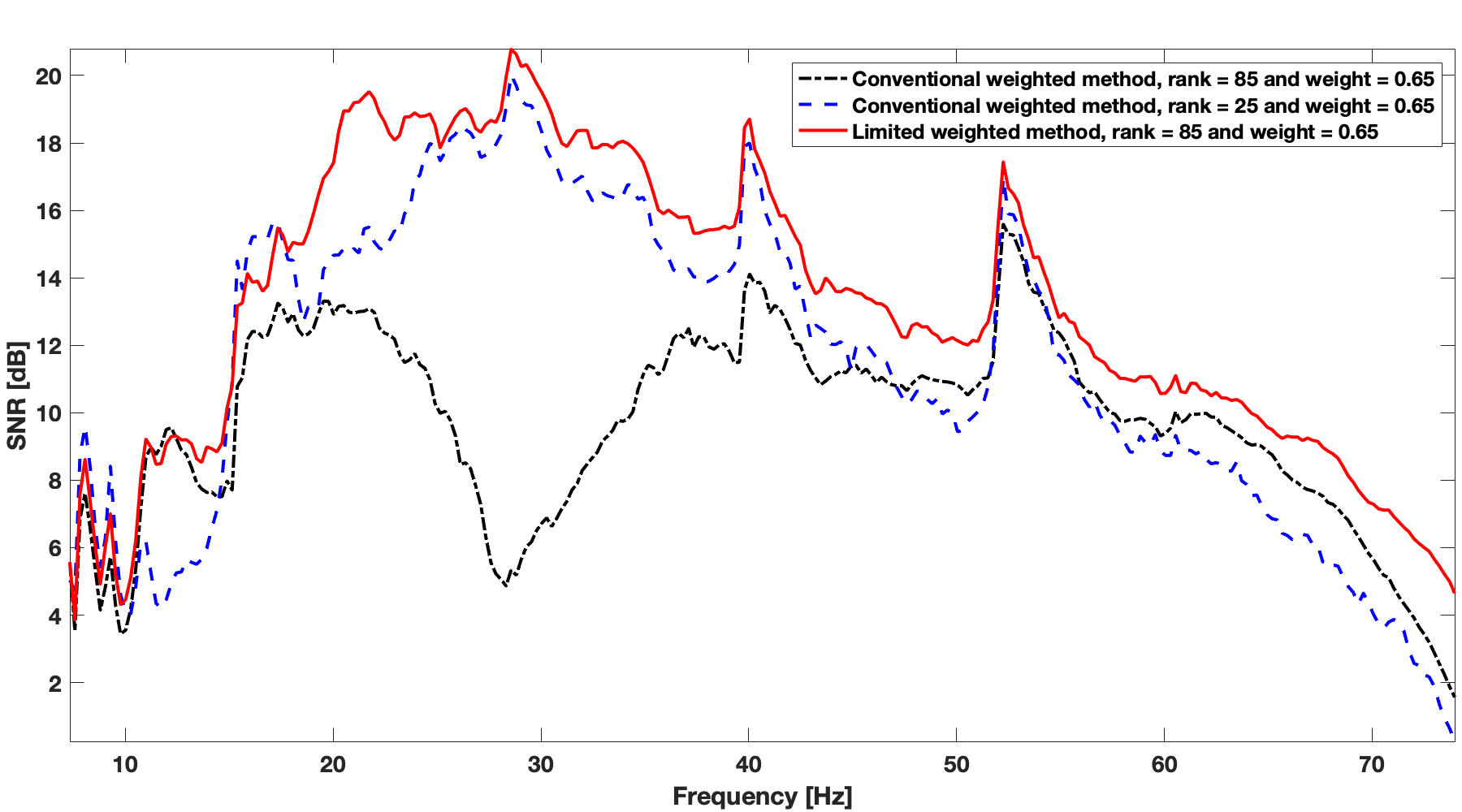}
\caption{$S/R$ of reconstructed data vs frequency based on our
limited-subspace weighted method (red color), conventional weighted
method with rank equals to $85$ (black color) and $25$ (blue
color).}\label{example_2}
\end{figure}

To further compare our limited-subspace method with the original method,
we repeat wavefield reconstructions over a range of frequencies
$7-74\, \mathrm{Hz}$. In Figure~\ref{example_2}, we show the comparison
of the S/R's across the whole frequency range. As expected, we observe
that limited-subspace weighted method (red line in
Figure~\ref{example_2}) outperforms conventional weighted method for
both ranks of $25$ (blue line in Figure~\ref{example_2}) and $85$ (black
line in Figure~\ref{example_2}) for most of the frequencies. This is
because of using limited subspace we avoid risk of overfitting at lower
frequencies and hence get improvement in quality of reconstructed data.

\begin{figure}
\centering
\subfloat[\label{full_time1}]{\includegraphics[width=0.500\hsize]{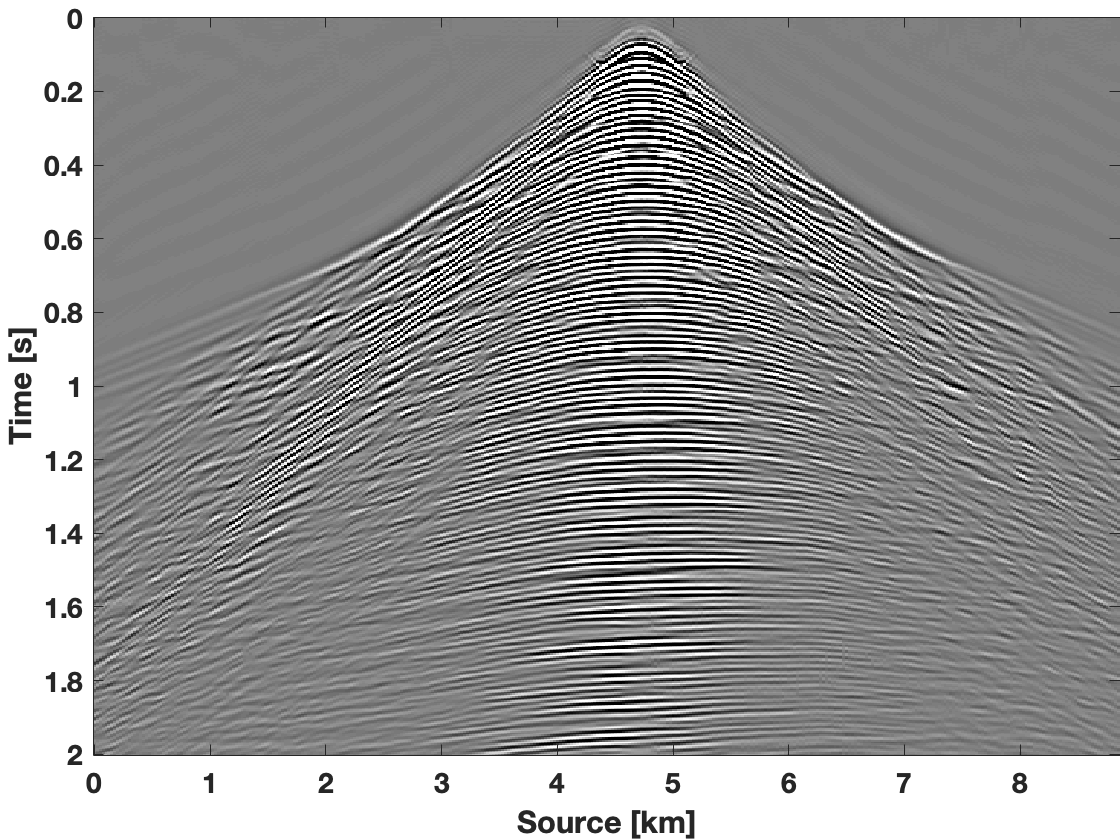}}
\subfloat[\label{jittered_time_75}]{\includegraphics[width=0.500\hsize]{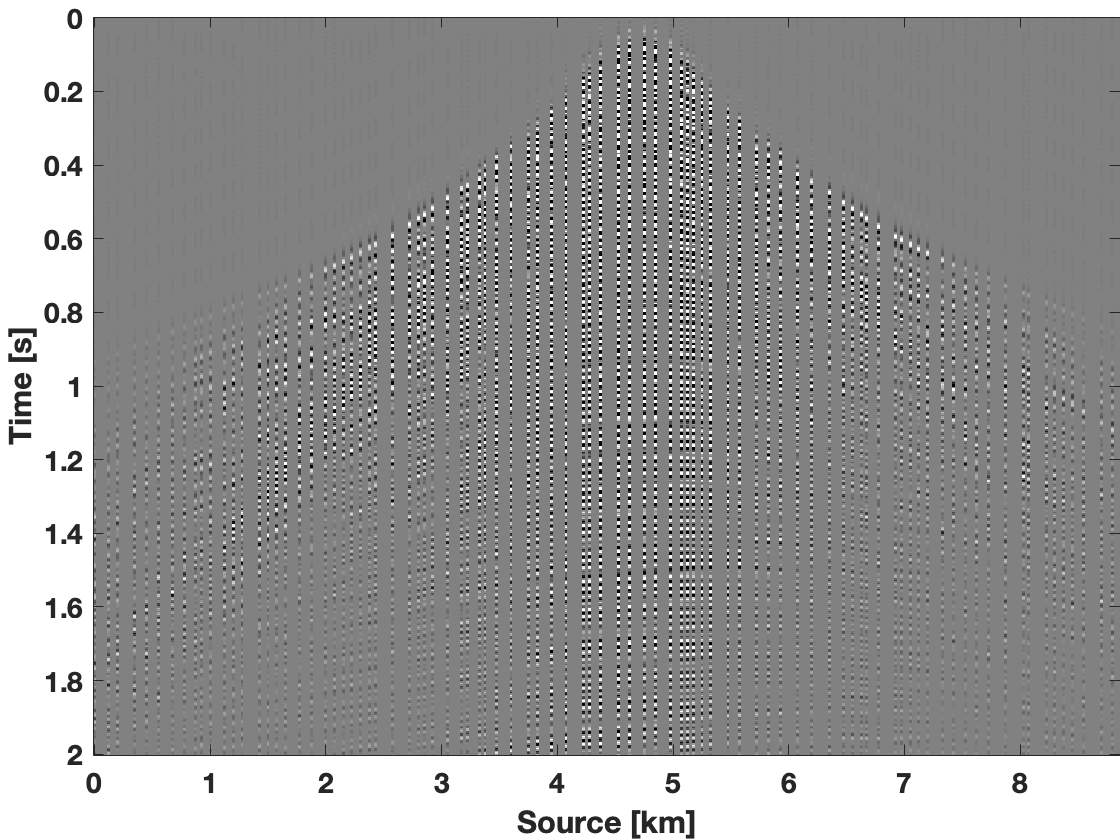}}
\\
\subfloat[\label{forward_time_diff}]{\includegraphics[width=0.500\hsize]{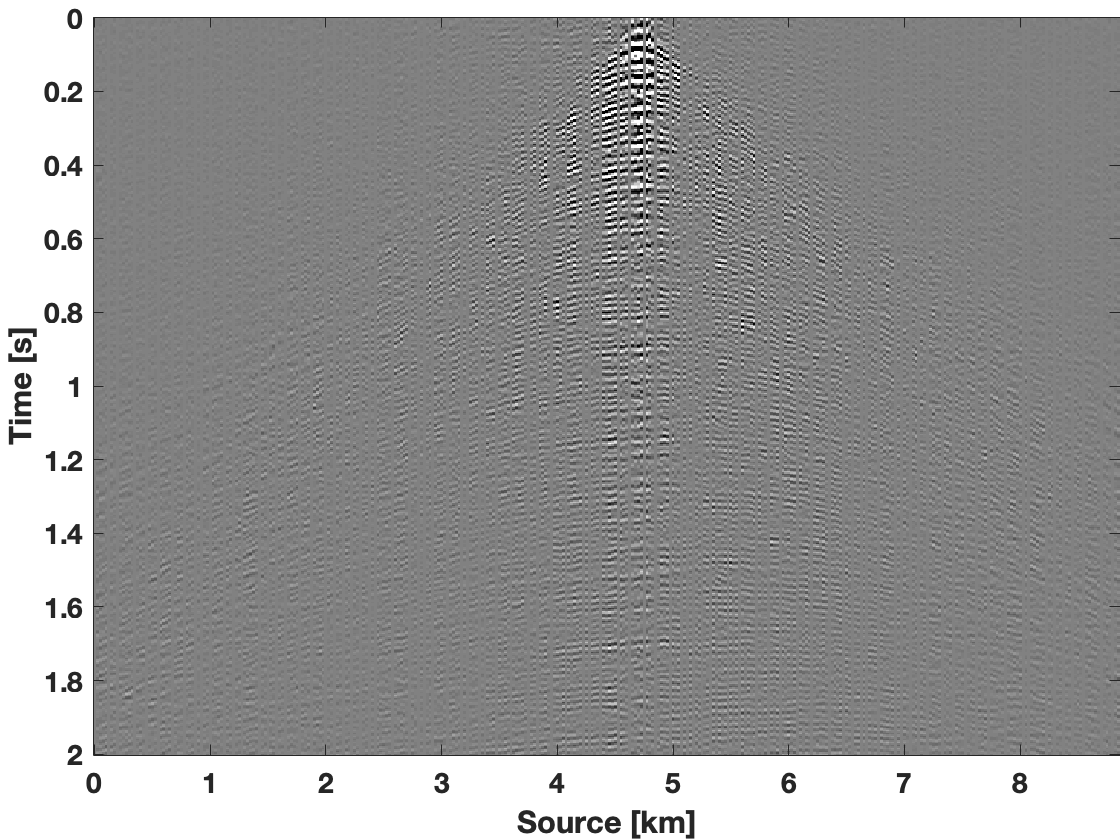}}
\subfloat[\label{reverse_time_diff}]{\includegraphics[width=0.500\hsize]{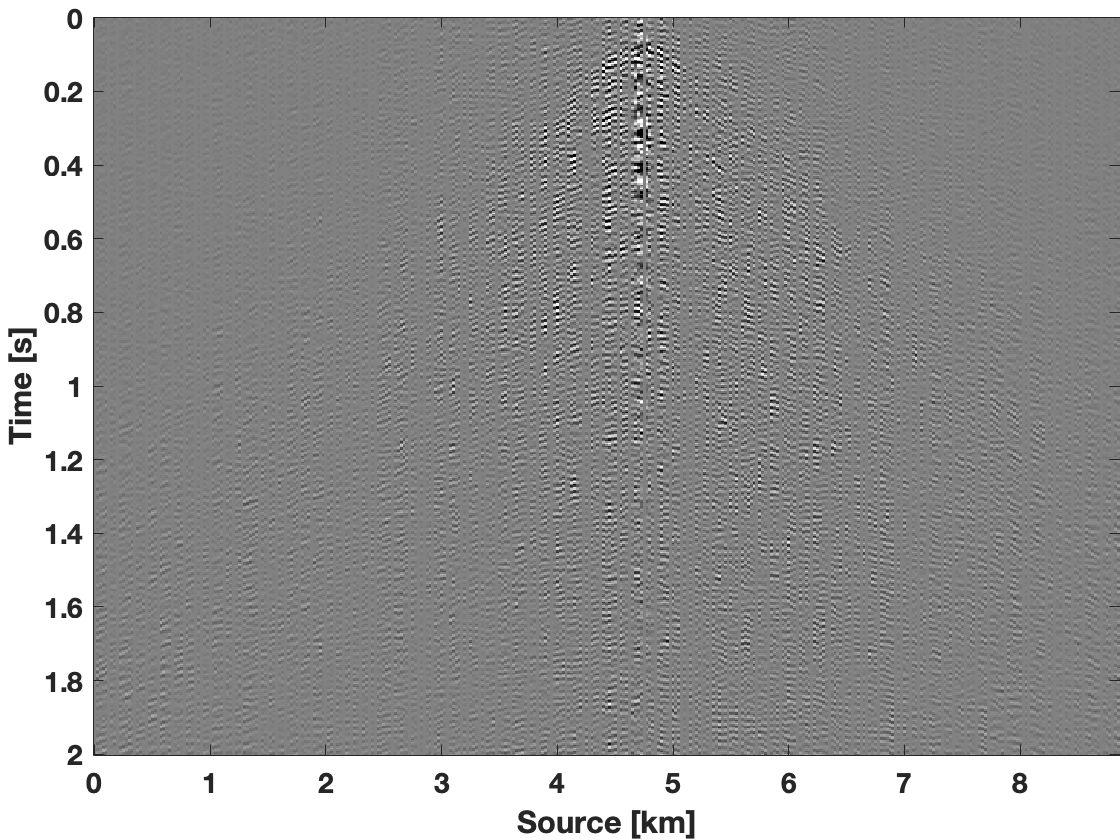}}
\\
\subfloat[\label{forward_time_85_diff}]{\includegraphics[width=0.500\hsize]{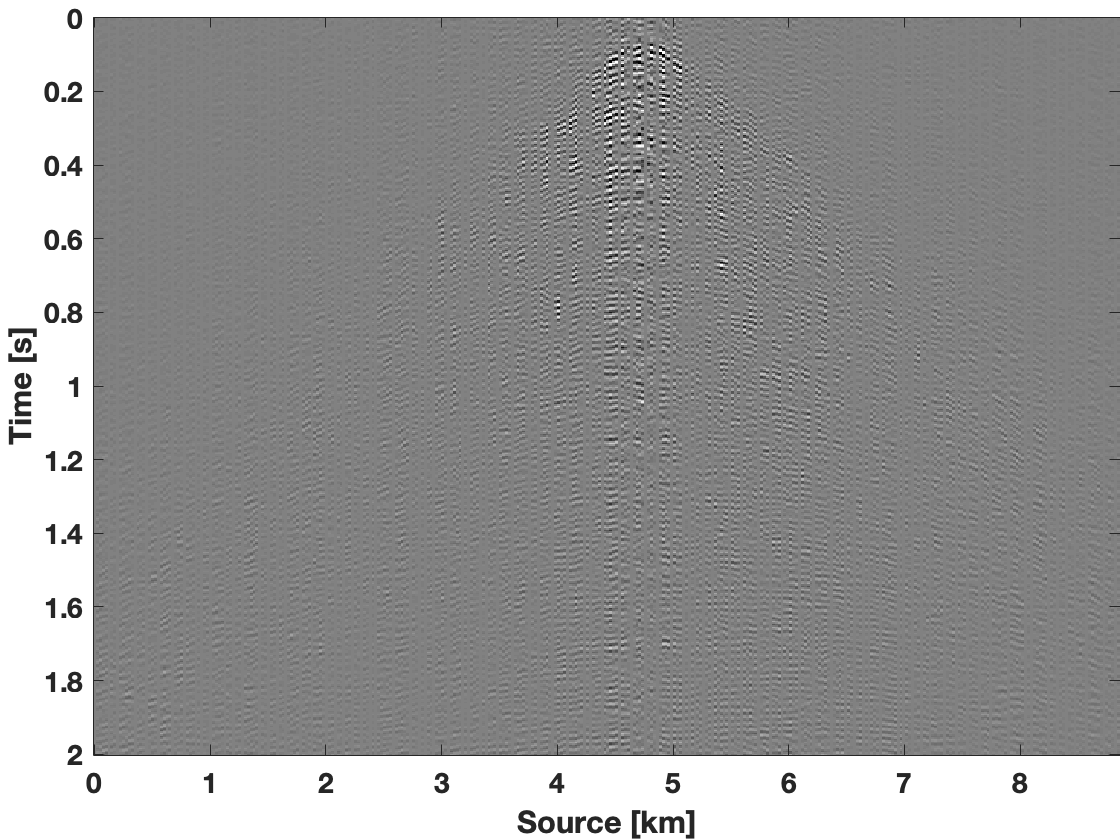}}
\caption{Wavefield reconstruction results in the time-domain. \emph{(a)}
Ground truth. (b) $75 \%$ subsampled seismic data with jittered
subsampling. \emph{(c)} using conventional weighted method
($S/R = 10.69\, \mathrm{dB}$) for rank equals to $r=85$, \emph{(d)}
using conventional weighted method ($S/R = 11.49\, \mathrm{dB}$) for
rank equals to $r=25$, \emph{(e)} using limited subspace weighted method
($S/R = 13.31\, \mathrm{dB}$) with limited subspace rank
$r_s=25$.}\label{example_4}
\end{figure}

To show the recovery improvement in the time domain, we included
Figure~\ref{example_4}. To make fair comparison, we construct a bandpass
filter with pass frequency $7-74\, \mathrm{Hz}$ with a transition width
at both ends of $3.66\, \mathrm{Hz}$. We apply this bandpass filter on
the true data, the subsampled data, and on recovered data recovered
using the three scenarios described above. After applying the filter, we
transform the filtered data back to the time domain. As we can see from
Figure~\ref{forward_time_85_diff}, we observe less leakage of coherent
signal in the data residual for results obtained with our
limited-subspace weighted method in comparison to the data residual
yielded by the conventional weighted method with ranks of $r=85$
(Figure~\ref{forward_time_diff}) and $r=25$
(Figure~\ref{reverse_time_diff}). With the conventional weighted method
for rank equals to $r=85$, we get S/R of $10.69\, \mathrm{dB}$, and for
rank $r=25$, we get S/R of $11.49\, \mathrm{dB}$. With the
limited-subspace weighted method we get S/R of $13.31\, \mathrm{dB}$,
which is a significant improvement.

\vspace*{-0.45cm}

\section{Conclusions}\label{conclusions}

In this work, we proposed a limited-subspace weighted method to further
improve the performance of recursively weighted method in terms of
better data reconstruction quality. By exploiting the fact that
dimensions of weight matrices are independent of the rank of the
subspaces, our method allows us to use higher ranks for data
reconstruction while avoiding the risk of overfitting at the lower
frequencies. Matrices with higher rank allow for a better approximation
of the frequency slices at higher frequencies and hence allow for better
quality of reconstructed data if we prevent overfitting by working with
limited-subspace weights. Through experiments we performed on a field
data acquired in the Gulf of Suez, we demonstrated the advantage of our
method in comparison to the recursively weighted method without using
limited subspace. We also introduced a computationally more efficient
formulation by moving the weight matrices to the data-misfit term. In
future work , we would like to extend the application of
limited-subspace weighted method to large scale $3$D data examples.

\vspace*{-0.45cm}

\section{Related materials}\label{related-materials}

In order to facilitate the reproducibility of the results herein
discussed, Matlab \& Julia implementation of this work are made
available on the \href{https://slim.gatech.edu}{SLIM} GitHub page
\url{https://github.com/slimgroup/Software.SEG2020}.

\vspace*{-0.45cm}

\section{Acknowledgement}\label{acknowledgement}

We would like to acknowledge the support from Georgia Institute of
Technology for funding this research.

\bibliography{abstract}

\end{document}